\documentclass[11pt]{article} 
\usepackage{epsfig} 
\usepackage{amssymb,amsmath,amsthm,amscd}
\usepackage{latexsym} 
\usepackage{cancel}
\usepackage{xcolor}
\usepackage{amssymb}

\usepackage[margin=.92in]{geometry}
\usepackage{subcaption}
\usepackage{graphicx}
\usepackage{url}
\usepackage{authblk}
 
\newtheorem{thm}{Theorem}[section] 
\newtheorem{prop}[thm]{Proposition}

\newtheorem{conj}{Conjecture}[section]   
 
\theoremstyle{definition}

\theoremstyle{remark} 
\newtheorem{rem}{Remark}[section]  
\def\eqref#1{(\ref{#1})} 

\newcommand{\bbR}{{\mathbb{R}}}

\newcommand {\mat}      [1] {\left[\begin{array}{#1}}
\newcommand {\rix}          {\end{array}\right]}
\newcommand {\de}      [1] {\left|\begin{array}{#1}}
\newcommand {\nt}          {\end{array}\right|}

\newcommand{\bstar}       {\begin{eqnarray*}}
\newcommand{\estar}       {\end{eqnarray*}}
\newcommand{\eqn}       {\begin{eqnarray}}
\newcommand{\enn}       {\end{eqnarray}}
\newcommand{\eq}[1]   {\begin{equation}\label{#1}}
\newcommand{\en}      {\end{equation}}

\begin{document}
\begin{titlepage}
\title{Bifurcation of Flux Ratio in Ionic Flows via a PNP model} 
\author[1,2]{Hamid Mofidi\thanks{\tt h.mofidi@bimsa.cn}}
\affil[1]{Beijing Institute of Mathematical Sciences and Applications (BIMSA), Beijing 101408, China}
\affil[2]{Yau Mathematical Sciences Center, Tsinghua University, Beijing 100084, China}

\date{}
\end{titlepage}

\maketitle

{\textbf{Abstract.}}
   \smallskip
This study investigates the flux ratios of ionic flows using Poisson-Nernst-Planck (PNP) models. The flux ratio measures the effect of permanent charge on the ion fluxes and can exhibit bifurcation, where a small parameter change causes a sudden change in system behavior. We identify the bifurcation points of flux ratios and analyze their dependence on boundary conditions and channel geometry. We use numerical simulations to verify our theoretical results and to explore the qualitative behaviors of fluxes and flux ratios at the bifurcation points. We provide insights into the dynamics of ion channels and their applications in ion channel design and optimization.
 \\


\smallskip


\vspace{5mm}

\section{Introduction.}\label{Intro}
\setcounter{equation}{0} 
Ion channels are membrane proteins that enable the diffusion of ions across cell membranes, generating electric signals for cellular communication. Their types and functions depend on their permanent charge distributions and channel shapes, which are the main structural features. The interaction of these features with physical parameters, such as boundary electric potentials and concentrations of ionic mixtures, determines the complex dynamics of ionic flows. Mathematical models based on the Poisson-Nernst-Planck (PNP) equations are essential for studying and analyzing these dynamics.
 This model incorporates the interaction of structural features with physical parameters and has been extensively analyzed using a geometric singular perturbation approach \cite{ Eis1, Eis4, EL07}.

A key measure in the study of ion channels is the flux ratio, which quantifies the relative movement of different ion species across a membrane channel. The flux ratio reflects the effect of permanent charge on the ion fluxes and can indicate the selectivity, permeability, and reversal potential of the channel. The flux ratio can also exhibit bifurcation, where a small parameter change causes a sudden change in system behavior.
The flux ratio is influenced by various factors, such as boundary concentrations, diffusion coefficients, ion size and valence, and different potentials that model steric and electrostatic interactions. These factors affect the driving force and the resistance for ion movement, resulting in complex effects on the flux ratio \cite{GB08, JEL19, NC18}.

The interplay of these factors results in complex effects on the flux ratio. Mathematical models based on the PNP equations, which describe the electrodiffusion of ions in a fluid, are crucial for studying and analyzing these factors. These models can be solved analytically or numerically to obtain expressions and values of the flux ratio under different conditions \cite{EL07,JEL19}.

\subsection{Main Findings and Their Sections.}
This manuscript investigates the flux ratios of ionic flows using PNP models. The flux ratio measures the relative movement of different ion species across a membrane channel and is influenced by the permanent charge, boundary conditions, and channel geometry. The flux ratio can also exhibit bifurcation, where a small parameter change causes a sudden change in system behavior. We combine theoretical and numerical approaches to identify and analyze the bifurcation points of flux ratios and their dependence on various factors. We also explore the qualitative behaviors of fluxes and flux ratios at the bifurcation points. We provide insights into the dynamics of ion channels and their applications in ion channel design and optimization.  The main results and their corresponding sections are summarized as follows:

Section \ref{sec-PNP} introduces the theoretical framework based on the PNP systems for ion channels. Section \ref{sec-setup} presents the basic setup and the relevant results that are used throughout the paper, including the mathematical model derived from the PNP equations, the GSP theory and its application to PNP, and the definition of the flux ratios. Section \ref{sec-Fluxbif} focuses on the bifurcation behavior of the flux ratios and examines the conditions and factors that affect the occurrence of bifurcation moments in ionic flows. Section \ref{sec-Numbifflux} extends the investigation through numerical simulations and provides a quantitative and qualitative analysis of flux ratios and bifurcation phenomena in ionic flows. Section \ref{sec-discussion}, Discussion and Future Work, explores potential avenues for future research and discusses the implications and applications of the findings in various fields.

	\subsection{ PNP Systems for Ion Channels.}\label{sec-PNP}
	\setcounter{equation}{0}
	

The PNP equations have been extensively simulated and computed in previous studies \cite{ CE,EL07}. These simulations have highlighted the need to include mathematical boundary conditions, such as macroscopic reservoirs, to accurately describe the behavior of ion channels \cite{GNE, NCE}.
In the PNP model for an ionic mixture of $n$ ion species, denoted by $k=1,2,\cdots, n$, the equations are as follows:

Poisson Equation:
\begin{align}
\nabla \cdot \left( \varepsilon_r(\overrightarrow{X}) \varepsilon_0 \nabla \Phi \right) = -e_0 \left( \sum_{s=1}^n z_s C_s + \mathcal{Q}(\overrightarrow{X}) \right),
\end{align}
where $\overrightarrow{X}$ is within a three-dimensional cylindrical-like domain $\Omega$ representing the channel, $\varepsilon_r(\overrightarrow{X})$ is the relative dielectric coefficient, $\varepsilon_0$ is the vacuum permittivity, $e_0$ is the elementary charge, and ${\mathcal{Q}}(\overrightarrow{X})$ represents the permanent charge density of the channel in molar units (M).

Nernst-Planck Equations:
\begin{align}
\partial_t C_k + \nabla \cdot \overrightarrow{\mathcal{J}}_k = 0, \quad -\overrightarrow{\mathcal{J}}_k = \dfrac{1}{k_BT} \mathcal{D}_k(\overrightarrow{X}) C_k \nabla \mu_k,
\end{align}
where $C_k$ is the concentration of the $k$-th ion species, $z_k$ is the valence, and $\mu_k$ is the electrochemical potential. The flux density $ \overrightarrow{{\mathcal{J}}}_k(\overrightarrow{X})$ represents the number of particles crossing each cross-section per unit time, ${\mathcal{D}}_k(\overrightarrow{X})$ is the diffusion coefficient, and $n$ is the number of distinct ion species.

	In the subsequent section, we will delve into the discussion of quasi-one-dimensional models, which are derived from three-dimensional PNP systems due to the relatively thin cross-sections of ion channels compared to their lengths \cite{LW10}.


\section{Basic Setup and Relevant Results.}\label{sec-setup}
\setcounter{equation}{0}
This section provides a comprehensive overview of our mathematical model for ionic flows, encompassing the fundamental setup and pertinent outcomes. Specifically, we investigate a quasi-one-dimensional PNP model that characterizes the ion transport within a confined channel featuring a permanent charge. 

To establish a solid foundation for our subsequent analysis, we introduce the notation and assumptions that will be consistently employed throughout the paper. Additionally, we incorporate a review of pertinent findings from previous literature that serve as vital underpinnings for our analysis. These previous works are important for understanding our main contributions, which we will present in the following sections  \cite{EL07, Liu05}.


\subsection{A Quasi-One-Dimensional PNP and a Dimensionless Form.}\label{sec-Quasi1}
Our analysis is based on a   quasi-one-dimensional PNP model first proposed in \cite{NE} and, for a special case, rigorously justified in \cite{LW10}.
For a mixture of  $n$ ion species,
a quasi-one-dimensional    PNP  model   is
\begin{align}\begin{split}\label{ssPNP}
&\frac{1}{A(X)}\frac{d}{dX}\Big(\varepsilon_r(X)\varepsilon_0A(X)\frac{d\Phi}{dX}\Big)=-e_0\Big(\sum_{s=1}^nz_sC_s+{\mathcal{Q}}(X)\Big),\\
 & \frac{d{\mathcal{J}}_k}{dX}=0, \quad -{\mathcal{J}}_k=\frac{1}{k_BT}{\mathcal{D}}_k(X)A(X)C_k\frac{d\mu_k}{d X}, \quad
 k=1,2,\cdots, n,
\end{split}
\end{align}
where $X\in [a_0,b_0]$ is the coordinate along the axis of the channel and baths of total length $b_0-a_0$, $A(X)$ is the
 area of cross-section  of the channel over the longitudinal location $X$, $e_0$ is the elementary charge (we reserve the letter $e$ for the Euler's  number -- the base for the natural exponential function), $\varepsilon_0$ is the vacuum permittivity, $\varepsilon_r(X)$ is the relative dielectric coefficient, ${\mathcal{Q}}(X)$ is the permanent charge density, $k_B$ is the Boltzmann constant, $T$ is the absolute temperature, $\Phi$ is the electric potential,   for the $k$th ion species, $C_k$ is the concentration, $z_k$ is the valence, ${\mathcal{D}}_k(X)$ is the diffusion coefficient, $\mu_k$ is the electrochemical potential, and ${\mathcal{J}}_k$ is the flux density.

Equipped with the system (\ref{ssPNP}),  a meaningful    boundary condition  for ionic flow through ion channels (as explained in \cite{EL07}) can be expressed as follows for $k=1,2,\cdots, n$
\begin{equation}
\Phi(a_0)={\mathcal{V}}, \ \ C_k(a_0)={\mathcal{L}}_k>0; \quad \Phi(b_0)=0,  \ \
C_k(b_0)={\mathcal{R}}_k>0. \label{ssBV}
\end{equation}

In relation to typical experimental designs, the positions $X=a_0$ and $X=b_0$ are located in the baths separated by the channel and are locations for two electrodes that are applied to control or drive the ionic flow through the ion channel. 
An important measurement is the {\em I-$\mathcal{V}$ (current-voltage) relation} where, for fixed ${\mathcal{L}}_k$'s and ${\mathcal{R}}_k$'s, the current $I$ depends on the transmembrane potential (voltage) ${\mathcal{V}}$:
\[
I=\sum_{s=1}^nz_s{\mathcal{J}}_s({\mathcal{V}}).
\] 
Of course, the relations of individual fluxes ${\mathcal{J}}_k$'s on ${\mathcal{V}}$ contain more information but it is much harder to experimentally measure the individual fluxes ${\mathcal{J}}_k$'s.
Ideally, the experimental designs should not affect the intrinsic ionic flow properties so one would like to design the boundary conditions to meet the so-called electroneutrality 
\[\sum_{s=1}^nz_s{\mathcal{L}}_s=0=\sum_{s=1}^nz_s{\mathcal{R}}_s.\]
The reason for this is that, otherwise, there will be sharp boundary layers which cause significant changes  (large gradients) of the electric potential and concentrations near the boundaries so that a measurement of these values has non-trivial uncertainties. 
One smart design to remedy this potential problem is the ``four-electrode-design": two `outer electrodes' in the baths far away from the ends of the ion channel to provide the driving force and two `inner electrodes'' in the bathes near the ends of the ion channel to measure the electric potential and the concentrations as the ``real'' boundary conditions for the ionic flow. At the inner electrodes locations, the electroneutrality conditions are reasonably satisfied, and hence, the electric potential and concentrations vary slowly and a measurement of these values would be robust.

The following rescaling or its variations have been widely used for the convenience of mathematical analysis \cite{Gil99}.   
  Let $C_0$ be a characteristic concentration of the ion solution. 
 We now make a dimensionless re-scaling of the variables in the system (\ref{ssPNP}) as follows.
\begin{align}\label{rescale}\begin{split}
&\varepsilon^2=\frac{\varepsilon_r\varepsilon_0k_BT}{e_0^2(b_0-a_0)^2C_0},\; x=\frac{X-a_0}{b_0-a_0},\;  h(x)=\frac{A(X)}{(b_0-a_0)^2},
\;  Q(x)=\frac{{\mathcal{Q}}(X)}{C_0},  \\
&D(x)={\mathcal{D}}(X),\; \phi(x)=\frac{e_0}{k_BT}\Phi(X), \; c_k(x)=\frac{C_k(X)}{C_0}, \;  
 J_k=\frac{{\mathcal{J}}_k}{(b_0-a_0)C_0   {\mathcal{D}}_k}. 
\end{split}
\end{align}
  We assume $C_0$ is fixed but large so that the parameter $\varepsilon$ is small. Note that $\varepsilon=\frac{\lambda_D}{b_0-a_0},$
   where $\lambda_D$  is the Debye screening length.
   In terms of the new variables,    the BVP (\ref{ssPNP}) and (\ref{ssBV}) becomes 
\begin{align}\label{PNP2}
\begin{split}
&\frac{\varepsilon^2}{h(x)}\frac{d}{dx}\left(h(x)\frac{d\phi}{dx}\right)=-\sum_{s=1}^nz_s
c_s -Q(x),\\
&\frac{d J_k}{dx}=0, \quad -  J_k=\frac{1}{k_BT}D(x)h(x)c_k\frac{d \mu_k}{d x},
\end{split}
\end{align}
with  boundary conditions at $x=0$ and $x=1$
\begin{align}\label{BVO}
\begin{split}
\phi(0)=&V,\; c_k(0)=L_k; \;
 \phi(1)=0,\; c_k(1)=R_k,
\end{split}
\end{align}
where 
\[V:=\frac{e_0}{k_BT}{\mathcal{V}},\quad L_k:=\frac{{\mathcal{L}}_k}{C_0},\quad R_k:=\frac{{\mathcal{R}}_k}{C_0}.\]
The permanent charge $Q(x)$ is
\begin{align}\label{Q}
Q(x)=\left\{\begin{array}{cc}
0, & x\in (0,a)\cup (b,1)\\
2Q_0, &x\in (a,b),
\end{array}\right.
\end{align}
where 
\[0<a=\frac{A-a_0}{a_1-a_0}<b=\frac{B-a_0}{a_1-a_0}<1.\]
We will take the ideal component $\mu_k^{id}$ only for the electrochemical potential. In terms of the new variables, it becomes
\[\frac{1}{k_BT}\mu_k^{id}(x)=z_k\phi(x)+\ln c_k(x).\]

In this study, we will explore the boundary value problem (BVP) \eqref{PNP2} and its boundary conditions \eqref{BVO} for ionic mixtures featuring a cation with valence $z_1=1$ and an anion with valence $z_{2}=-1$. 
Our critical assumption is that $\varepsilon$ is small. This assumption allows us to treat the BVP (\ref{PNP2}) with (\ref{BVO}) as a singularly perturbed problem. A general framework for analyzing such singularly perturbed BVPs in PNP-type systems has been developed in prior works \cite{EL07,JLZ15,Liu05,LX15} for classical PNP systems and in \cite{JL12,LLYZ13,LM22,LTZ12} for PNP systems with finite ion sizes.

\subsection{An Overview of a GSP and Governing System for BVP (\ref{PNP2}) and (\ref{BVO}).} \label{RelResult}
 
To establish the foundation for our work, we begin by revisiting the geometric singular perturbation framework introduced in \cite{EL07}. Specifically, in the case where $n=2$ with $z_1>0>z_2$, the authors of \cite{EL07} employed geometric singular perturbation theory to construct the singular orbit for the boundary value problem (BVP) consisting of (\ref{PNP2}) and (\ref{BVO}).
Introducing the transformations $u=\varepsilon \phi$ and $w=x$, we can reformulate system (\ref{PNP2}) as follows for $k=1,\cdots,n$,
 \begin{align}\label{PNPO}
 \begin{split}
 \varepsilon  \dot\phi=&u, \quad \varepsilon \dot u=-\sum_{s=1}^nz_sc_s-Q(w)-\frac{\varepsilon h'(w)}{h(w)}u,\\
 \varepsilon\dot{c}_k=&-z_kc_ku -   \frac{\varepsilon J_k}{D(w)h(w)},\quad 
 \dot J_k=0,\quad \dot w=1,
\end{split}
\end{align}
where the symbol ``dot" denotes differentiation with respect to the variable $x$. The autonomous system (\ref{PNPO}) can then be treated as a dynamical system with phase space of $\bbR^{2n+3}$ and state variables $(\phi, u, c_1,\cdots, c_n, J_1, \cdots, J_n, w)$.
In view of the jumps of permanent charge $Q(x)$ at $x=a$ and $x=b$, the construction of singular orbits is split into three intervals $[0,a]$, $[a,b]$, $[b,1]$ as follows. To do so, one introduces (unknown) values of $(\phi,c_1,c_2)$ at $x=a$ and $x=b$:
\begin{align}\label{6unknown}
\phi(a)=\phi^a,\; c_1(a)=c_1^a,\; c_2(a)=c_2^a;\quad \phi(b)=\phi^b,\; c_1(b)=c_1^b,\; c_2(a)=c_2^b.
\end{align}
These values determine 
(boundary) conditions at $x=a$ and $x=b$  as 
\[B_a=\{(\phi^a, u, c_1^a, c_2^a, J_1, J_2, a):\;\mbox{arbitrary } u, J_1,J_2\},\]
and 
\[ B_b=\{(\phi^b, u, c_1^b, c_2^b, J_1, J_2, b): \;\mbox{arbitrary } u, J_1,J_2\}.\]
  This setup leads to six unknowns $\phi^{a}$, $\phi^b$, $c_k^a$ and $c_k^b$ for $k=1,2$ should be determined.
On each interval, a singular orbit typically consists two singular layers and one regular layer. 
\begin{description}
\item{(1)}\ On interval $[0,a]$, a singular orbit from $B_0$ to $B_a$ consists of two singular layers located at $x=0$ and $x=a$, denoted as $\Gamma_0^l$ and $\Gamma_a^l$, and one regular layer $\Lambda_l$. Furthermore, with the preassigned values $\phi^a$, $c_1^a$ and $c_2^a$, the flux $J_k^l$ and $u_l(a)$ are uniquely determined so that
\[(\phi^a,u_l(a), c_1^a, c_2^a, J_1^l, J_2^l, a)\in B_a.\]

\item{(2)}\ On interval $[a,b]$, a singular orbit from $B_a$ to $B_b$ consists of two singular layers located at $x=a$ and $x=b$, denoted as $\Gamma_a^r$ and $\Gamma_b^l$, and one regular layer $\Lambda_m$. Furthermore, with the preassigned values $(\phi^a, c_1^a, c_2^a)$ and $(\phi^b, c_1^b, c_2^b, b)$, the flux $J_k^m$, $u_m(a)$ and $u_m(b)$ are uniquely determined so that
\[(\phi^a,u_m(a), c_1^a, c_2^a, J_1^m, J_2^m, a)\in B_a\;\mbox{ and }\; (\phi^b,u_m(b), c_1^b, c_2^b, J_1^m, J_2^m, b)\in B_b.\]
\item{(3)}\ On interval $[b,1]$, a singular orbit from $B_b$ to $B_1$ consists of two singular layers are located at $x=b$ and $x=1$, denoted as $\Gamma_b^r$ and $\Gamma_1^l$, and one regular layer $\Lambda_r$. Furthermore, with the preassigned values $\phi^b$, $c_1^b$ and $c_2^b$, the flux $J_k^r$ and $u_r(b)$ are uniquely determined so that
\[(\phi^b,u_r(b), c_1^b, c_2^b, J_1^r, J_2^r, b)\in B_b.\]
\end{description}
 The matching conditions of this connecting problem are 
 \begin{align}\label{GSys}
 J_k^l=J_k^m=J_k^r\;\mbox{ for }\; k=1,2,\;   u_l(a)=u_m(a)\;\mbox{ and }\; u_m(b)=u_r(b).
 \end{align}
  There are total six conditions, which are exactly the same number of unknowns preassigned in (\ref{6unknown}). Then the singular connecting problem is reduced to {\em the governing system} (\ref{GSys}) (see \cite{EL07} for an explicit form of the governing system).


\subsection{Flux Ratios for Permanent Charge Effects.}\label{FR4Q}
The main goal in studying ion channels is to understand how the flows of different ion species depend on the channel's structure and boundary conditions. Experimental data typically provides only the total current, making it challenging to determine the individual ion flows.
To address this challenge, mathematical models, like the PNP model, relate individual ion flows to factors such as channel structure, permanent charge density, diffusion coefficients, and concentration gradients. Solving these models helps estimate ion flows under various conditions and gain insights into ion channel dynamics and the effects of permanent charge on each ion species \cite{HLY21, SL18}.

In \cite{Liu17},  to characterize the effects of permanent charges on fluxes for given boundary conditions, the  author introduces a ratio
\begin{equation}\label{lambda}
\lambda_k(Q;\varepsilon)=\frac{J_k(Q;\varepsilon)}{J_k(0;\varepsilon)},
\end{equation}
where $J_k(Q;\varepsilon)$ is the flux of $k$-th ion species associated to the permanent charge $Q(x)$ and $J_k(0;\varepsilon)$ is the flux associated to zero permanent charge. Since permanent charges cannot change the sign of flux, one has $\lambda_k(Q;\varepsilon)>0$.
If $\lambda_k(Q;\varepsilon)>1$, then the permanent charge $Q$  enhances the flux in the sense that $|J_k(Q;\varepsilon)|>|J_k(0;\varepsilon)|$.   On the other hand, if  $\lambda_k(Q;\varepsilon)<1$, then the permanent charge $Q$  reduces the flux in the sense that $|J_k(Q;\varepsilon)|<|J_k(0;\varepsilon)|$. 

In \cite{Liu17},  for mixtures of two ion species with $z_1>0>z_2$,   the following  universality of a permanent charge  effect   is established as follows. For ionic flow with one cation and one anion, under some general conditions but independent of boundary conditions,  
\begin{align}\label{universal}
\mbox{\it if }\; Q(x)\ge 0,\;\mbox{\it then, for }\; \varepsilon>0 \;\mbox{\it small, }
\lambda_1(Q;\varepsilon)\le \lambda_2(Q;\varepsilon),
\end{align}
where $\lambda_1(Q;\varepsilon)$ is the ratio associated to the cation and $\lambda_2(Q;\varepsilon)$ is the ratio associated to the anion.
Furthermore, the statement (\ref{universal}) is sharp in the sense that, depending on the boundary conditions,  each one of the followings is possible (\cite{JLZ15})
\begin{itemize}
\item[(i)] $1<\lambda_1(Q;\varepsilon)<\lambda_2(Q;\varepsilon)$ (both cation and anion fluxes are enhanced); 
\item[(ii)] $\lambda_1(Q;\varepsilon)<1<\lambda_2(Q;\varepsilon)$ (cation flux is reduced but anion flux is enhanced);
\item[(iii)] $\lambda_1(Q;\varepsilon)< \lambda_2(Q;\varepsilon)<1$ (both cation and anion fluxes are reduced).
\end{itemize}

\section{Bifurcations of Critical Flux Ratios.}\label{sec-Fluxbif}
\setcounter{equation}{0}

  For fixed $L$ and $R$, the critical ratio is $\lambda_k(V,Q)=1$.  We will be interested in the bifurcation of $\lambda_k(V,Q)=1$; that is, we are interested in the bifurcation of the solution $V=V_k(Q)$ of $\lambda_k(V,Q)=1$ viewing $Q$ as a parameter.   
     It follows that the necessary condition for  bifurcations at $(V^*,Q^*)$ are: 
\begin{equation}\label{bif.cond.}
\lambda_k(V^*,Q^*)=1\;\mbox{ and }\;  \partial_V\lambda_k(V^*,Q^*)=0.
\end{equation}
Bifurcations are phenomena where the qualitative behavior of a dynamical system changes as a parameter varies. There are different types of bifurcations, depending on other conditions, but they all require to satisfy \eqref{bif.cond.} for some fixed point $V^*$ and parameter value $Q^*$, where $\lambda_k(V, Q)$ is the flux ratio function. 
\begin{rem}\label{biftypes}
\em
We focus on identifying general bifurcation moments in ionic flows. Various types of bifurcations can occur, but we won't delve into specific details here. Our main aim is to find the bifurcation moments $(V^*, Q^*)$ and analyze the relationships between variables at these points.
For example, a saddle-node bifurcation occurs when stable and unstable fixed points are created or destroyed as $Q$ crosses a critical value. Similarly, a pitchfork bifurcation results in three equilibrium points, with specific conditions. Transcritical and Hopf bifurcations involve their own stability exchanges and periodic orbits.
While examining all these conditions simultaneously can be complex, we will not explore these intricacies in this paper due to their computational demands.
\end{rem}

\subsection{Preparation of Bifurcation Moments for $n=2$ with $z_1=1=-z_2$.}\label{sec-prep}

We now consider $\lambda_k(Q_0)$ for $n=2$ with $z_1=-z_2=1$.  For simplicity, we assume    electroneutrality boundary conditions
$L_1=L_2=L$ and $R_1=R_2=R$ in the following. 

In \cite{EL07}, for $z_1=1$ and $z_2=-1$, the governing system (\ref{GSys})  is reduced to an equation with only one unknown $A$. More precisely, 
set $c_1^{a}c_{2}^{a}=A^2$ and $c_1^{b}c_{2}^{b}=B^2$ with $A>0$, $B>0$,
\[\alpha=\frac{H(a)}{H(1)}\;\mbox{ and }\; \beta=\frac{H(b)}{H(1)}\;\mbox{ where }\; H(x)=\int_0^x\frac{1}{D(s)h(s)}ds.\]
Introduce $I=J_1-J_2$ and $F=J_1+J_2$ so that, for $k=1,2$,
\[J_k=\frac{1}{2}\left(F+(-1)^{k+1}I\right).\]
The governing system becomes
\begin{align}\label{0eF} 
  e^{-Fy}\Big(F\sqrt{Q_0^2+A^2}+Q_0I\Big) -F\sqrt{Q_0^2+B^2}-Q_0I=0,
\end{align}
where $B$, $y$, $I$ and $F$ are determined, for $Q_0\neq 0$, in terms of the variable $A$ by
\begin{align}\label{0ev}
\begin{split}
B=&\frac{1-\beta}{\alpha}(L-A)+R,\quad F=2\frac{L-A}{H(a)},\\
I=&-\frac{2(L-A)}{H(a)\ln{\frac{BL}{AR}} }\left(\ln{\frac{B}{A}}-V-\ln\frac{\sqrt{Q_0^2+B^2}-Q_0}{\sqrt{Q_0^2+A^2}-Q_0}\right.\\
&\left.\quad -\frac{\sqrt{Q_0^2+B^2}-\sqrt{Q_0^2+A^2}}{Q_0}-\frac{(\beta-\alpha)(L-A)}{\alpha Q_0}\right),\\
Iy=&-\frac{(\beta-\alpha)(L-A)}{\alpha Q_0}+\frac{\sqrt{Q_0^2+A^2}-\sqrt{Q_0^2+B^2}}{Q_0}.
 \end{split}
\end{align}
 
\noindent In the scenario where $Q_0=0$ and $z_1=-z_2=1$, the expressions for $J_k(0)$ with $k=1,2$ can be found in Proposition 3.1 of the reference \cite{JLZ15}. Under electroneutrality conditions, these expressions are,

\begin{align}\label{zeroQJ}
J_k(0) = \frac{(L-R)\left((-1)^{k+1}V + \ln L - \ln R\right)}{H(1)(\ln L - \ln R)}.
\end{align}

\noindent Therefore, the values of $\lambda_k$ in equation \eqref{lambda} for $k=1,2$ can be expressed as

\begin{equation}\label{lambda2}
\lambda_k(Q_0) = \frac{J_k(Q_0)}{J_k(0)} = \frac{\frac{1}{2}\left(\frac{2(L-A)}{H(a)} + (-1)^{k+1}I\right)}{(L-R)\left((-1)^{k+1}V + \ln L - \ln R\right)/(H(1)(\ln L - \ln R))}.
\end{equation}

Recall from \cite{LX15} that if $Q_0=0$, then the BVP has a unique solution. We thus consider only the case where $Q_0\neq 0$. 
Set
\begin{equation}\label{Q-free}
{\mathcal{A}}:=\frac{A}{Q_0},\; {\mathcal{B}}:=\frac{B}{Q_0};\; l:=\frac{L}{Q_0},\;r:=\frac{R}{Q_0},\; f:=\frac{H(a)F}{2Q_0}, \; {\mathcal{I}}:=\frac{H(a)I}{2Q_0},\; Y:=\frac{2Q_0y}{H(a)}.
\end{equation}

\begin{rem}\label{signs}
  \em Note that $\mathcal{A}$, $\mathcal{B}$, $l$, $r$ and $Y$   should have the same sign as that of $Q_0$, and $Q_0$ could be negative. 
\end{rem}

\begin{rem}\label{Positive Q}
\em We assume, without loss of generality, that  $Q_0$ is positive. In other words, the PNP for two ion species has the symmetry with respect to $Q_0$ to $- Q_0$, $z_1$ to $- z_1$, $z_2$ to $- z_2$, and $\phi$ to $- \phi$.
\end{rem}
In terms of the new variables introduced in \eqref{Q-free}, the governing system \eqref{0eF} will be,
\begin{align}\label{1eF}
 {\mathcal{I}}\ln\frac{{\mathcal{I}}-({\mathcal{A}}-l)\sqrt{1+{\mathcal{B}}^2}}{{\mathcal{I}}-({\mathcal{A}}-l)\sqrt{1+{\mathcal{A}}^2}}
=\frac{\beta-\alpha}{\alpha}({\mathcal{A}}-l)^2+\big(\sqrt{1+{\mathcal{A}}^2}-\sqrt{1+{\mathcal{B}}^2}\big)({\mathcal{A}}-l),
\end{align}
together with
\begin{align}\label{t2i}
{\mathcal{I}}=&\frac{{\mathcal{A}}-l}{\ln{\frac{{\mathcal{B}}l}{{\mathcal{A}}r}} }\left(\ln{\frac{{\mathcal{B}}}{{\mathcal{A}}}}-V-\ln\frac{\sqrt{1+{\mathcal{B}}^2}-1}{\sqrt{1+{\mathcal{A}}^2}-1} - \sqrt{1+{\mathcal{B}}^2}+\sqrt{1+{\mathcal{A}}^2}+\frac{(\beta-\alpha)({\mathcal{A}}-l)}{\alpha}\right).
\end{align}
 For the moment, we fix $l$, $r$, $\alpha$ and $\beta$. From \ref{t2i}, we can write ${\mathcal{I}}={\mathcal{I}}({\mathcal{A}},V)$, and then, from \eqref{1eF}, we expect to solve ${\mathcal{A}}={\mathcal{A}}(V)$.
 For this purpose, we set
\begin{align}\label{2eF}
{\mathcal{F}}({\mathcal{A}},{\mathcal{I}}):= {\mathcal{I}}\ln\frac{{\mathcal{I}}-({\mathcal{A}}-l)\sqrt{1+{\mathcal{B}}^2}}{{\mathcal{I}}-({\mathcal{A}}-l)\sqrt{1+{\mathcal{A}}^2}}
-\frac{\beta-\alpha}{\alpha}({\mathcal{A}}-l)^2-\big(\sqrt{1+{\mathcal{A}}^2}-\sqrt{1+{\mathcal{B}}^2}\big)({\mathcal{A}}-l).
\end{align}
Note that ${\mathcal{F}}$ is defined for ${\mathcal{I}}$ that satisfies
\begin{equation}\label{I-bound}
{\mathcal{I}}>({\mathcal{A}}-l)\cdot\max\{\sqrt{1+{\mathcal{A}}^2},\sqrt{1+{\mathcal{B}}^2}\} \;\mbox{ or }\;  {\mathcal{I}}<({\mathcal{A}}-l)\cdot\min\{\sqrt{1+{\mathcal{A}}^2},\sqrt{1+{\mathcal{B}}^2}\}.
\end{equation}
Additionally, ${\mathcal{F}}$ is defined for ${\mathcal{A}}$ in the range $[0,{\mathcal{A}}_M]$, where ${\mathcal{A}}_M= \dfrac{\alpha r}{1-\beta}+l$, ensuring that ${\mathcal{B}} \geq 0$.

\subsection{Bifurcation Moments for PNP Models.}\label{sec-bifanalysis}
In this section, we study the bifurcation moments of the relation $\lambda_k(V)=1$ where $\lambda_k(V)$ is the flux ratio defined in \eqref{lambda}. We recall that $\lambda_k(V)$ depends on the variables ${\mathcal{A}}$ and ${\mathcal{I}}$, which are solutions of the governing system  \eqref{1eF} and \eqref{t2i}. We also recall that $k=1$ or $2$ corresponds to the cation or anion species, respectively.
For $k=1,2$,  we introduce the following notation from equation \eqref{lambda2},
\begin{align}\label{j_k}
j_k(V) &= \frac{H(a)}{2Q_0}J_k({\mathcal{A}}(V)) \nonumber \\
&= \frac{1}{2}\left(l - {\mathcal{A}} + (-1)^{k+1}{\mathcal{I}}\right).
\end{align}
where ${\mathcal{A}}={\mathcal{A}}(V)$ is a solution of the governing system  \eqref{1eF} and ${\mathcal{I}}={\mathcal{I}}({\mathcal{A}}(V),V)$ is given by \eqref{t2i}. Note that, from \eqref{zeroQJ}, we have
 \begin{equation}\label{j_0}
 \begin{aligned}
 j_k(0):=&\frac{H(a)}{2Q_0}J_k(0)=\frac{\alpha(l-r)((-1)^{k+1}V+\ln l-\ln r)}{2(\ln l-\ln r)},\\
 \frac{\partial j_k(0)}{\partial V}=&(-1)^{k+1}\frac{\alpha(l-r)}{2(\ln l-\ln r)}.
 \end{aligned}
  \end{equation}
 Thus,
 \begin{equation}\label{lambda3}
 \begin{aligned}
 \lambda_k(V)=&\frac{j_k({\mathcal{A}}(V))}{j_k(0)},\\
 \partial_V\lambda_k(V)=&\frac{\frac{\partial j_k}{\partial {\mathcal{A}}}\cdot\frac{\partial {\mathcal{A}}}{\partial V}+\frac{\partial j_k}{\partial V}}{j_k(0)}-\frac{j_k({\mathcal{A}}(V))}{j_k^2(0)}\cdot\frac{\partial j_k(0)}{\partial V}.
 \end{aligned}
  \end{equation}

We say that $V^*$ is a {\em bifurcation moment} for the relation  $\lambda_k(V)=1$ (or equivalently $\bar{\lambda}_k(V)= \lambda_k(V) - 1 = 0$ ) if $\lambda_k(V^*)=1$ and the equation cannot be solved for $V$ uniquely near $V^*$.  In the sequel, we will denote $j_k=j_k({\mathcal{A}}(V))$. 
 Therefore, for $k=1$ or $2$, $V=V^*$ is a bifurcation moment if and only if $\lambda_k(V)=1$ and $\partial_V\lambda_k(V)=0$; that is,
\begin{equation}\label{bifeqn0}
\begin{aligned} 
&j_k({\mathcal{A}}(V))=j_k(0),\\
& \frac{\partial j_k}{\partial {\mathcal{A}}} \frac{\partial {\mathcal{A}}}{\partial V}+\frac{\partial j_k}{\partial {\mathcal{I}}}(\frac{\partial {\mathcal{I}}}{\partial {\mathcal{A}}}\frac{\partial {\mathcal{A}}}{\partial V}+\frac{\partial {\mathcal{I}}}{\partial V})  = \frac{\partial j_k(0)}{\partial V}.
\end{aligned}
\end{equation}
To solve the system given by equation \eqref{bifeqn0}, we differentiate ${\mathcal{F}}({\mathcal{A}},{\mathcal{I}})=0$ with respect to $V$, where ${\mathcal{F}}$ is defined in equation (\ref{2eF}). 
Then by conducting a meticulous computation, for $k=1$ or $2$, the bifurcation moment is determined by the values of $({\mathcal{A}}, {\mathcal{I}},V)$ that satisfy the following system of equations:
\begin{equation}
\begin{aligned}\label{BifPara}
& {\mathcal{I}}\ln\frac{{\mathcal{I}}-({\mathcal{A}}-l)\sqrt{1+{\mathcal{B}}^2}}{{\mathcal{I}}-({\mathcal{A}}-l)\sqrt{1+{\mathcal{A}}^2}}
=\frac{\beta-\alpha}{\alpha}({\mathcal{A}}-l)^2+\big(\sqrt{1+{\mathcal{A}}^2}-\sqrt{1+{\mathcal{B}}^2}\big)({\mathcal{A}}-l),\\
&{\mathcal{I}}=\frac{{\mathcal{A}}-l}{\ln{\frac{ {\mathcal{B}}l}{{\mathcal{A}}r}} }\left(\ln\frac{{\mathcal{B}}\big(\sqrt{1+{\mathcal{A}}^2}-1\big)}{{\mathcal{A}}\big(\sqrt{1+{\mathcal{B}}^2}-1\big)} -V+ \sqrt{1+{\mathcal{A}}^2}-\sqrt{1+{\mathcal{B}}^2}+\frac{(\beta-\alpha)}{\alpha}({\mathcal{A}}-l)\right),\\
&  (-1)^{k+1}{\mathcal{I}} +l-{\mathcal{A}} =\frac{\alpha(l-r)((-1)^{k+1}V+\ln l-\ln r)}{\ln l-\ln r},\\
&\frac{\partial {\mathcal{F}}}{\partial {\mathcal{A}}}\Big(\frac{l-{\mathcal{A}}}{\ln\frac{l{\mathcal{B}}}{r{\mathcal{A}}}}-  \frac{\alpha(l-r)}{\ln l-\ln r}\Big)   =\frac{\partial \mathcal{F}}{\partial {\mathcal{I}}}\Big( (-1)^{k} \frac{l-{\mathcal{A}}}{\ln\frac{l{\mathcal{B}}}{r{\mathcal{A}}}}+ \frac{\alpha(l-r)}{\ln l-\ln r}\frac{\partial {\mathcal{I}}}{\partial {\mathcal{A}}}\Big),
\end{aligned}
\end{equation}
where ${\mathcal{F}}$ is given by \eqref{2eF}.

{ Note that the system (\ref{BifPara}) with four equations is overdetermined for three variables $({\mathcal{A}},{\mathcal{I} },V)$ as expected for bifurcation moment. It simply says that the system could be consistent only for special values of $Q_0$ encoded in $l$ and $r$. Due to this consideration, we can take $\sigma:=l/r=L/R$ as a fixed parameter and take $l$ as a free variable. Then the system (\ref{BifPara}) will become, for $(l,{\mathcal{A}},{\mathcal{I}},V)$,
\begin{equation}\label{2Bifreduced}
\begin{aligned}
& \rho = {\mathcal{I}}\ln\frac{{\mathcal{I}}-({\mathcal{A}}-l)\sqrt{1+{\mathcal{B}}^2}}{{\mathcal{I}}-({\mathcal{A}}-l)\sqrt{1+{\mathcal{A}}^2}}, \quad V=\ln\frac{{\mathcal{B}}\big(\sqrt{1+{\mathcal{A}}^2}-1\big)}{{\mathcal{A}}\big(\sqrt{1+{\mathcal{B}}^2}-1\big)} - \dfrac{\big({\mathcal{I}}\ln{\frac{\sigma {\mathcal{B}}}{{\mathcal{A}}}} - \rho\big)}{{\mathcal{A}}-l},\\
&  {\mathcal{I}} =\dfrac{\frac{(-1)^{k+1} \sigma\ln \sigma}{\alpha l(\sigma-1)} ({\mathcal{A}}-l)^2 + \big(\ln\frac{{\mathcal{B}}\big(\sqrt{1+{\mathcal{A}}^2}-1\big)}{{\mathcal{A}}\big(\sqrt{1+{\mathcal{B}}^2}-1\big)} - (-1)^{k} \ln \sigma \big)({\mathcal{A}}-l)+ \rho}{\ln{\frac{\sigma {\mathcal{B}}}{{\mathcal{A}}}} +  \frac{\sigma\ln \sigma}{\alpha l(\sigma-1)} ({\mathcal{A}}-l)},\\
&\Big( \frac{\partial \mathcal{F}}{\partial {\mathcal{I}}}-(-1)^{k}\frac{\partial {\mathcal{F}}}{\partial {\mathcal{A}}}\Big) V =\ln\frac{\sigma {\mathcal{B}}}{{\mathcal{A}}}\Big(\frac{\partial {\mathcal{F}}}{\partial {\mathcal{A}}} + \frac{\partial \mathcal{F}}{\partial {\mathcal{I}}} \frac{\partial {\mathcal{I}}}{\partial {\mathcal{A}}} \Big)- \ln \sigma\frac{\partial \mathcal{F}}{\partial {\mathcal{A}}}\\
&\hspace*{1.55in}  + (-1)^{k}\Big( {\mathcal{I}} \frac{\ln\frac{\sigma {\mathcal{B}}}{{\mathcal{A}}}}{{\mathcal{A}}-l}\big(\frac{\partial \mathcal{F}}{\partial {\mathcal{A}}} + \frac{\partial \mathcal{F}}{\partial {\mathcal{I}}} \frac{\partial {\mathcal{I}}}{\partial {\mathcal{A}}} \big) +  \ln \sigma\frac{\partial \mathcal{F}}{\partial {\mathcal{I}}}\Big),
\end{aligned}
\end{equation}
where, 
\begin{equation}\label{rho}
\rho= \rho ({\mathcal{A}},l) := \frac{\beta-\alpha}{\alpha}({\mathcal{A}}-l)^2+\big(\sqrt{1+{\mathcal{A}}^2}-\sqrt{1+{\mathcal{B}}^2}\big)({\mathcal{A}}-l),
\end{equation}
and $\sigma = \dfrac{l}{r}$, and $\rho$ is defined as in \eqref{rho}.
Set,
\begin{equation}\label{gamma}
\begin{aligned}
\gamma_1({\mathcal{I}},{\mathcal{A}},l):= \frac{1}{{\mathcal{I}}-({\mathcal{A}}-l)\sqrt{1+{\mathcal{A}}^2}}, \hspace*{.3in} \gamma_2({\mathcal{I}},{\mathcal{A}},l):= \frac{1}{{\mathcal{I}}-({\mathcal{A}}-l)\sqrt{1+{\mathcal{B}}^2}}.
\end{aligned}
\end{equation}
We then need to determine the expressions for $\frac{\partial {\mathcal{I}}}{\partial {\mathcal{A}}}$, $\frac{\partial \mathcal{F}}{\partial {\mathcal{A}}}$, and $\frac{\partial \mathcal{F}}{\partial {\mathcal{I}}}$ through careful computations. 
Thus, from \eqref{2Bifreduced} and the computations above, we have the following Proposition that provides the bifurcation moment through a system of algebraic equations as follows:

\begin{prop}\label{C4bifur}
The bifurcation moment for the flux ratio $\lambda_k,~k=1,2$, is determined by the values of $(l,{\mathcal{A}}, {\mathcal{I}},V)$ that satisfy the following system of equations,
\begin{equation}\label{1Bif-equiv}
\begin{aligned}
 \rho =&  {\mathcal{I}}\ln\frac{{\mathcal{I}}-({\mathcal{A}}-l)\sqrt{1+{\mathcal{B}}^2}}{{\mathcal{I}}-({\mathcal{A}}-l)\sqrt{1+{\mathcal{A}}^2}},\quad  V=\ln\frac{{\mathcal{B}}\big(\sqrt{1+{\mathcal{A}}^2}-1\big)}{{\mathcal{A}}\big(\sqrt{1+{\mathcal{B}}^2}-1\big)} - \dfrac{\big({\mathcal{I}}\ln{\frac{\sigma {\mathcal{B}}}{{\mathcal{A}}}} - \rho\big)}{{\mathcal{A}}-l},\\
 {\mathcal{I}} =&\dfrac{\frac{(-1)^{k+1} \sigma\ln \sigma}{\alpha l(\sigma-1)} ({\mathcal{A}}-l)^2 + \big(\ln\frac{{\mathcal{B}}\big(\sqrt{1+{\mathcal{A}}^2}-1\big)}{{\mathcal{A}}\big(\sqrt{1+{\mathcal{B}}^2}-1\big)} - (-1)^{k} \ln \sigma \big)({\mathcal{A}}-l)+ \rho}{\ln{\frac{\sigma {\mathcal{B}}}{{\mathcal{A}}}} +  \frac{\sigma\ln \sigma}{\alpha l(\sigma-1)} ({\mathcal{A}}-l)},\\
  \Big( &\frac{\beta-\alpha}{\alpha} - ({\mathcal{A}}-l)\big({\mathcal{A}}\gamma_1 + \frac{1-\beta}{\alpha}{\mathcal{B}}\gamma_2 \big) \Big) \Big(\ln\frac{\sigma {\mathcal{B}}}{{\mathcal{A}}}  +({\mathcal{A}}-l) \dfrac{\sigma\ln \sigma}{\alpha l(\sigma-1)}\Big)\\
  & \qquad \quad -\Big( \frac{\beta-\alpha}{\alpha} - ({\mathcal{A}}-l)\big({\mathcal{A}}\gamma_1+ \frac{1-\beta}{\alpha}{\mathcal{B}}\gamma_2 \big) \Big)\cdot M \\
  & \qquad  \quad +\dfrac{1}{{\mathcal{A}}-l}  \Big({\mathcal{I}} +  (-1)^k({\mathcal{A}}-l)\Big)\dfrac{\sigma\ln \sigma}{\alpha l(\sigma-1)}\cdot M\\
 & \qquad    =\frac{{\mathcal{I}}^2-({\mathcal{A}}-l)^2}{({\mathcal{A}}-l)}\cdot M\cdot \Big( \frac{1}{{\mathcal{A}}}\gamma_1 + \frac{(1-\beta)}{\alpha } \frac{1}{{\mathcal{B}}}\gamma_2\Big),
\end{aligned}
\end{equation}
where $M({\mathcal{I}},{\mathcal{A}},l) = {\mathcal{I}}\big(\gamma_2({\mathcal{I}},{\mathcal{A}},l)-\gamma_1({\mathcal{I}},{\mathcal{A}},l)\big) + \rho/{\mathcal{I}}$, and $\rho$ and $\gamma_i$ for $i=1,2$ were previously defined in \eqref{rho} and \eqref{gamma}, respectively. 
\end{prop}

\noindent System \eqref{1Bif-equiv} is a complex nonlinear system and making analytical solutions is too complicated or even impossible.  Numerical techniques can overcome this challenge and provide accurate approximations of solutions. Therefore, in the subsequent section, we employ numerical methods to find the solutions.  These methods allow us to explore the behavior of the system and gain valuable insights into the bifurcation moment.

\section{Numerical Study of Flux Ratio Bifurcation.}\label{sec-Numbifflux}
\setcounter{equation}{0}
In this section, we conduct numerical studies on the qualitative behaviors of flux ratios $\lambda_k$, which depend on several parameters (besides $Q_0$, including $(V, L, R, h)$. Our numerical investigation complements the analytical findings in \cite{JEL19,Liu17}, particularly in understanding flux ratios for small and large permanent charges. Analytical analysis is limited in providing even qualitative results for moderate permanent charge sizes. This section explores how fluxes and flux ratios change with different variables at the bifurcation point $V^*$ for a given $Q^*$ and under specific boundary conditions. Additionally, we present numerical results for $\lambda_1$ and $\lambda_2$ with fixed boundary concentrations of species.

In the previous section, we showed that for $k=1$ or $2$, the bifurcation moments are determined by the values of $(l,{\mathcal{A}}, {\mathcal{I}},V)$ that satisfy system \eqref{1Bif-equiv}, which we plan to solve numerically in this section.  The system \eqref{1Bif-equiv} is solved using the root function in Python from the scipy.optimize module, which finds the roots of a system of nonlinear equations. We use the `hybr' method, which stands for `hybrid' and refers to a solver that combines a modified Powell's method with a dogleg trust-region method. This method is suitable for systems of equations where the Jacobian matrix is either not available or approximated numerically. It is a good choice for general-purpose root-finding problems \cite{NW06}.
We use the following initial guesses for the unknown variables to solve the system of equations \eqref{1Bif-equiv}:
$
0<l\le 10, \quad  0<\mathcal{A} \le 10, \quad -60 \le \mathcal{I}\le 60,  \quad -80 \le V\le 80.$
However, these initial values are not fixed but depend on the concentrations obtained from \eqref{I-bound} and the constraint for $\mathcal{A}$ afterward, which vary with the parameter $r$. Therefore, the solutions and the ranges of the variables may change as we vary $r$. The figures show the numerical results of this variation.

\begin{rem} \em Note that we have the option to work with either the quantities $(A,B,L,R,I,J_k)$ or $(\mathcal{A}, \mathcal{B}, \mathcal{L}, \mathcal{R}, \mathcal{I}, j_k)$ as defined in \eqref{Q-free} and \eqref{j_k}. However, for the purpose of studying the qualitative behavior of the variables and for the sake of simplification, we will continue using the latter set. This means that given a specific $Q_0$, we are seeking the values of the other variables at the bifurcation moment. For accurate values, one can utilize equations \eqref{Q-free} and \eqref{j_k}.
\end{rem}

\subsection{Understanding Flux Ratios and Fluxes at the Bifurcation Points.}\label{sec-Num1}
Assuming $k=1$ in the system \eqref{1Bif-equiv}, we can determine the bifurcation points of the flux ratio $\lambda_1$ when $\lambda_1 =1$ remains constant, while $\lambda_2 > 1$ varies in relation to $V$ and $\mathcal{I}$, with $r$ fixed at $2$ and $0<l\leq 10$. The behavior of $\lambda_2$ is illustrated in Figure \ref{fig-lamVI}, where it is observed that $\lambda_2$ initially increases for negative values of $V$ and $\mathcal{I}$, reaching its maximum values at $(V^*,\mathcal{I}^*)$ for $V^*<0$ and $\mathcal{I}^*<0$. Subsequently, $\lambda_2$ decreases as $V$ and $\mathcal{I}$ increase towards zero and eventually assumes positive values. 
In this context, we propose the following conjecture, drawing from the insights provided by Figure \ref{fig-lamVI}:


\begin{conj}\label{conj-I}
At the bifurcation moment when the flux ratio $\lambda_1=1$, the system \eqref{1Bif-equiv} exhibits the following properties for the flux ratio as a function of $\mathcal{I}$ and $V$:

\begin{enumerate}
\item[i.] For a fixed boundary condition $r_0>0$, there exists a critical solution denoted as $(l^*, \mathcal{A}^*, V^*, \mathcal{I}^*)$ where $V_0^* < 0$ and $\mathcal{I}^* < 0$. Furthermore, the partial derivatives exhibit the following characteristics:

\begin{itemize}
\item[a.] $\frac{\partial \lambda_2}{\partial V} > 0$, ~ $\frac{\partial^2 \lambda_2}{\partial V^2} > 0$ \quad  for $V < V^*,$
\quad b. $\frac{\partial \lambda_2}{\partial \mathcal{I}} > 0$, ~ $\frac{\partial^2 \lambda_2}{\partial \mathcal{I}^2} > 0$ \quad for $\mathcal{I} <  \mathcal{I}^*$,
\item[c.] $\frac{\partial \lambda_2}{\partial V} < 0$, ~ $\frac{\partial^2 \lambda_2}{\partial V^2} > 0$ \quad  for $V > V^*,$
\quad d. $\frac{\partial \lambda_2}{\partial \mathcal{I}} < 0$, ~ $\frac{\partial^2 \lambda_2}{\partial \mathcal{I}^2} > 0$ \quad for $\mathcal{I} >  \mathcal{I}^*$.
\end{itemize}

\item[ii.] For fixed $0 < r_1 < r_2$, the corresponding critical solutions $(l_1^*, \mathcal{A}_1^*, V_1^*, \mathcal{I}_1^*)$ and $(l_2^*, \mathcal{A}_2^*, V_2^*, \mathcal{I}_2^*)$ satisfy $V_2^* < V_1^* < 0$ and $\mathcal{I}_2^* < \mathcal{I}_1^* < 0$. Furthermore, it holds that  $\lambda_2(V_2^*) < \lambda_2(V_1^*)$ and $\lambda_2(\mathcal{I}_2^*) < \lambda_2(\mathcal{I}_1^*)$.


\end{enumerate}
\end{conj}

\begin{figure}[htbp]
 	\centering
 	\begin{subfigure}{1.0\textwidth}
 		\includegraphics[width=\textwidth]{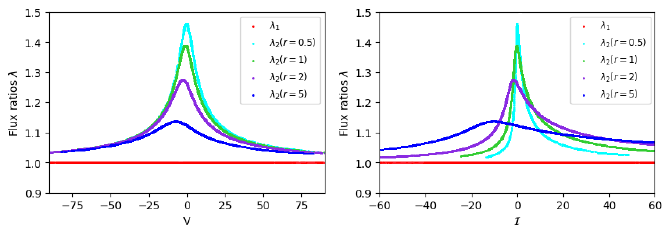}
 	\end{subfigure}	
 	\caption{\em The variation of the flux ratio $\lambda_2$ with respect to $V$ and $\mathcal{I}$ at the bifurcation moment for $\lambda_1=1$(horizontal red lines), for different values of the boundary concentration $r$. 
 	 }\label{fig-lamVI}
 \end{figure}

%

In this part, we investigate the influence of boundary concentrations $l$ and $r$ on  fluxes $j_k$'s, for $k=1,2$, and flux ratio $\lambda_2$ at the bifurcation moment of the other flux ratio, i.e., $\lambda_1$. Recall from \eqref{j_k} that for $k=1,2$,
$
j_k = \left(l - {\mathcal{A}} + (-1)^{k+1}{\mathcal{I}}\right)/2.
$
It is also noteworthy that in our calculations, we treated $l$ as a variable within the system described by equation \eqref{1Bif-equiv}, while $r$ was held constant. Figure \ref{fig-lflx} (left panel) illustrates critical values of $l$ where $\lambda_2$ undergoes a change in direction. Motivated by this observation, we present the following conjecture:


\begin{conj}\label{conj-L}
At the bifurcation moment when the flux ratio $\lambda_1 = 1$, the solutions of the system \eqref{1Bif-equiv} exhibit the following characteristics for the boundary concentrations:

\begin{enumerate}
\item[i.] For a fixed boundary condition $r_0 > 0$, there exists a critical solution $(l^*, \mathcal{A}^*, V^*, \mathcal{I}^*)$ where $l^* > 0$. Moreover, the partial derivatives satisfy the following:

\begin{itemize}
\item[a.] $\frac{\partial \lambda_2}{\partial l} > 0$, $\frac{\partial^2 \lambda_2}{\partial l^2} < 0$ \quad for $0 < l < l^*$,
\qquad b. $\frac{\partial \lambda_2}{\partial l} < 0$, $\frac{\partial^2 \lambda_2}{\partial l^2} > 0$ \quad for $l > l^*$,
\item[c.] $j_1 < 0$, $j_2 > 0$ \quad for $0 < l < l^*$, \qquad \quad  ~~d. $j_1 > 0$, $j_2 < 0$ \quad for $l > l^*$,
\item[e.] The flux $j_1$ increases for any boundary concentration $l$, while $j_2$ decreases, regardless of the value of $l$. In other words, $\frac{\partial j_1}{\partial l} > 0$ and $\frac{\partial j_2}{\partial l} < 0$ for any $l > 0$.
\end{itemize}

\item[ii.] For fixed $0 < r_1 < r_2$, the corresponding critical solutions $(l_1^*, \mathcal{A}_1^*, V_1^*, \mathcal{I}_1^*)$ and $(l_2^*, \mathcal{A}_2^*, V_2^*, \mathcal{I}_2^*)$ satisfy $0 < l_2^* < l_1^*$. Furthermore, $\lambda_2(l_2^*) < \lambda_2(l_1^*)$. However, $l^*$ is the critical point of the fluxes where $j_n(l^*) = 0$ for $n=1,2$.
\end{enumerate}
\end{conj}

\begin{figure}[htbp]
 	\centering
 	\begin{subfigure}{1.0\textwidth}
 		\includegraphics[width=\textwidth]{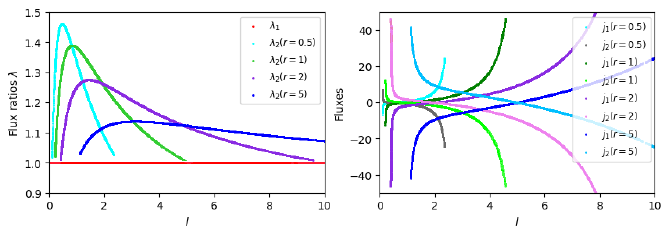}
 	\end{subfigure}	
 	\caption{\em The variation of the flux ratio $\lambda_2$ and fluxes $j_1$ and $j_2$ as a function of boundary concentration $l$ at the bifurcation moment for $\lambda_1=1$, for different values of the boundary concentration $r$. The figure illustrates that $\lambda_2$ reaches its maximum values at $l^*>0$. 
 }\label{fig-lflx}
 \end{figure}

\subsection{Understanding the Interplay of Electric Potential, Current Effects, and Boundary Concentration at $\lambda_1$ Bifurcations }\label{sec-Num2}
In this section, we delve into the dynamic relationship between electric potential, current, and fluxes at $\lambda_1$ bifurcations. Our aim is to understand how these variables interplay and impact the system's behavior. 
It is important to note that the current, denoted as $\mathcal{I}$, is not a constant but rather a variable that depends on several other parameters. Specifically, we explore the relationship $\mathcal{I} = z_1j_1 + z_2j_2$, which, under the condition $z_1=-z_2=1$, simplifies to $\mathcal{I} = j_1 - j_2$ is in alignment with the insights presented in Figure \ref{fig-flxIV}.

Figure \ref{fig-flxIV} visually represents the dynamics of these relationships, with distinct patterns emerging. We observe that 
$
\frac{\partial j_1}{\partial \mathcal{I}}> 0, \quad  \frac{\partial j_2}{\partial \mathcal{I}}< 0,\quad  \frac{\partial j_1}{\partial V}> 0, \quad  \frac{\partial j_2}{\partial V}< 0.
$
These observations form the basis for our exploration in this section, as we seek to uncover the underlying mechanisms governing electric potential and current effects on fluxes at $\lambda_1$ bifurcations.  

\begin{figure}[htbp]
 	\centering
 	\begin{subfigure}{1.0\textwidth}
 		\includegraphics[width=\textwidth]{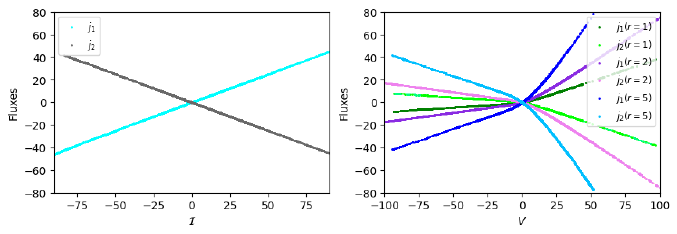}
 	\end{subfigure}	
 	\caption{\em Flux dependence on the current $\mathcal{I}$ and tranmembrane potential $V$ at $\lambda_1$ bifurcation moment for various values of $r$. }\label{fig-flxIV}
 \end{figure}

Our numerical findings demonstrate that the current $\mathcal{I}$ and the fluxes $j_k$ remain unaffected by the boundary concentrations $l$ and $r$. This consistency with the theoretical relationship $\mathcal{I} = j_1 - j_2$ is expected by definition.
Figure \ref{fig-flxIV} displays a single plot of $\mathcal{I}$ versus $j_k,$ as the other plots exhibit a similar pattern, with the only variation being in the range of $\mathcal{I}$, which depends on the value of $r.$ Finally, we explore how the current $\mathcal{I}$ and boundary electric potential $V$ behave in relation to the boundary concentration $l,$ examining different scenarios with varying values for the second boundary concentration $r.$

\begin{figure}[htbp]
 	\centering
 	\begin{subfigure}{1.0\textwidth}
 		\includegraphics[width=\textwidth]{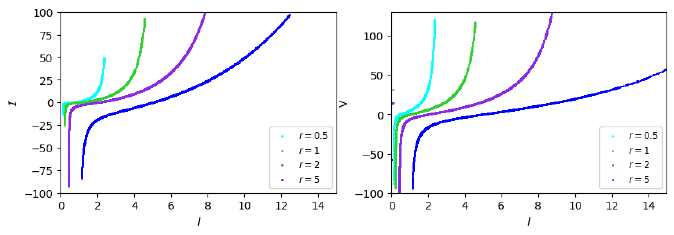}
 	\end{subfigure}	
 	\caption{\em Current $I$ and boundary electric potential $V$ behavior for various values of the boundary concentration $l$  at the bifurcation moment. }\label{fig-lIV}
 \end{figure}

 Figure \ref{fig-lIV} illustrates that both $V$ and $\mathcal{I}$ consistently share the same sign, irrespective of the boundary concentrations $l$ and $r$. 
This result matches well with the results of previous studies \cite{ELX15}, supporting our findings.


\section{Conclusions and Further Research.}\label{sec-discussion}

In this manuscript, we have investigated flux ratios in ionic flows with permanent charge, providing valuable insights into how it influences ion fluxes. We have uncovered a bifurcation phenomenon, where small parameter changes lead to significant system alterations. Employing PNP models, we have explored the bifurcation points and their dependencies on factors like boundary conditions and channel geometry. Our numerical results confirm that, for positive permanent charges, cation flux ratios are consistently smaller than anion flux ratios at the bifurcation point. We have combined theoretical analysis with numerical methods to gain deeper insights into ion channel dynamics and their practical applications \cite{FLMZ22, ML19, MEL20}.

Our work enhances the understanding of ionic flow mechanisms and their applications in channel design. Future research directions include exploring more complex models considering ion size effects, hard-sphere electrochemical potentials, or various permanent charge profiles. Additionally, we aim to investigate the impact of diffusion coefficients on bifurcation moments and delve into zero-current cases. For these numerical investigations, we'll employ diverse methods like root solvers, optimization algorithms, and neural networks to provide comprehensive insights into channel behavior and performance across different scenarios \cite{JLZ15, Liu17, M21, M20, ML19}.



\bibliographystyle{plain}

\newpage

\end{document}